%4ctJR.tex: 
%%a Plain TeX file by Shalosh B. Ekhad and Doron Zeilberger (2 page)
 
%begin macros

\baselineskip=14pt
\parskip=10pt

\magnification=\magstephalf

\def\1{{\overline{1}}}
\def\2{{\overline{2}}}
\parindent=0pt
\overfullrule=0in

\def\frac#1#2{{#1 \over #2}}
%\headline={\rm  \ifodd\pageno  \RightHead  \else  \LeftHead  \fi}
%\def\RightHead{\centerline{
%Title
%}}
%\def\LeftHead{ \centerline{Doron Zeilberger}}
%end macros
\centerline
{\bf 
Solving  Functional Equations Dear to  W.T. Tutte using the Naive 
}
\centerline
{\bf (yet fully rigorous!) Guess And Check Method}
\bigskip
\centerline
{\it Shalosh B. EKHAD and Doron ZEILBERGER}

{\bf Abstract}: In his seminal paper ``A census of planar triangulations", published in 1962, the iconic graph theorist (and code-breaker),
W.T. Tutte, spent a few pages to prove that a certain bi-variate generating function that enumerates triangulations, satisfies
a certain functional equation. He then used his genius to actually solve it, giving closed-form solutions to the
enumerating sequences. While the first part, of deriving the functional equation,  still needs human ingenuity, the
second part, of solving it, can nowadays be fully automated. Our Maple program, accompanying this paper, Tutte.txt,
can not only solve Tutte's original equation in a fraction of a second, it can also solve many, far more complicated ones,
way beyond the scope of even such a giant as W.T. Tutte. We use our favorite method of ``guess and check"
and show how it can always be made fully rigorous (if desired).

{\bf W.T. Tutte's functional equation}

In his seminal paper [T], W.T. Tutte introduced a certain bivariate generating function $\psi=\psi(x,y)$ whose coefficient of $x^n\, y^m$, $\psi_{n,m}$,
(in the Taylor expansion around $(0,0)$) enumerates planar triangulations with $3n+m$ internal edges and $m+3$ external edges.
The most interesting case was with $m=0$, i.e. when the perimeter is a triangle.

He first used his genius to derive a {\bf functional equation} (Eq. (3.8) of [T], p. 27), where $\psi=\psi(x,y)$ and $g=\psi(x,0)$:
$$
y^2 \psi^2 + (x+xgy-y-y^2)\psi + y - xg \, = \, 0.
$$
He then went on, using more than two pages, and lots of previous knowledge, to prove the expliict expression for $g(x)$:
$$
g\,=\, 1+ x + 2 \sum_{n=2}^{\infty} \frac{x^n}{(n+1)!} (3n+3)(3n+4) \cdots (4n+1) \quad .
$$

{\bf  Our Way}

Using our Maple package, {\tt Tutte.txt}, available from 

{\tt https://sites.math.rutgers.edu/\~{}zeilberg/tokhniot/Tutte.txt}  \quad,

and typing

{\tt Paper0(y**2*psi**2+(x+x*g*y-y-y**2)*psi+y-x*g,psi,g,x,y,5,n,1000);}

rederives, in less than a second, Tutte's result that indeed
$$
\psi_{n,0} \,= \,2\frac{ (3n+3)(3n+4) \cdots (4n+1)}{(n+1)!} \quad ,\quad ( n \geq 1) \quad .
$$

See the output file:

{\tt https://sites.math.rutgers.edu/\~{}zeilberg/tokhniot/oTutte1.txt}  \quad .

{\bf Our General Approach}

Let's describe our {\it naive} (but ingenious!) method, that is closely related to the method in [GZ].

{\bf Input}:  a polynomial in four variables $Q$.

{\bf Output}: a linear recurrence equation with polynomial coefficients for the coefficients $\psi_{n,0}$, of the unique bi-variate formal power series
$\psi(x,y)$ satisfying
$$
Q(\psi(x,y),\psi(x,0),x,y) \, = \, 0 \quad .
\eqno(1)
$$

Note that we assume that the functional equation  is well-defined, i.e. it uniquely defines a formal power series in $x,y$. This can be
easily checked by the computer.

{\bf Step 1}: Decide on a ``guessing parameter", {\tt K}, a positive integer, and use the functional equation (iteratively) to find the first $K$ coefficients, in $x$, of $\phi(x,y)$.
Note that the coefficient of each $x^i$ is a certain {\it rational function} in $y$.

{\bf Step 2}: Plug in $y=0$, getting the first $K$ coefficients (now these are numbers) of $g(x)=\psi(x,0)$.

{\bf Step 3}: Using the Maple command {\tt gfun[listoalgeq]} (see [SZ]) (or our own implementation in this package called {\tt Empir(L,x,P)}), {\it guess} an algebraic equation of the form
$$
P_1 (g(x),x) \, = \, 0 \quad .
\eqno(2)
$$

If $K$ is too small, you would get {\tt FAIL}. Then don't despair, just make $K$ larger.

Note that, so far it is only a guess! How do we prove it? (rigorously!).

Assuming that this equation is correct, we combine it with the functional equation $Q(\psi(x,y), g(x),x,y)=0$ and {\tt eliminate} $g(x)$, getting a 
{\bf polynomial} equation

$$
P_2(\psi(x,y),x,y)=0 \quad .
\eqno(3)
$$

We now {\bf claim} that the unique solution of $(3)$ satisfies $(1)$. But this is (rigorously) automatically provable. One way is to stay in
the {\bf algebraic ansatz}  and argue that $(3)$ (that implies $(1)$) implies a certain (complicated!) polynomial equation 
$P_4(Q,x,y)=0$, satisfied by
$Q(\psi,g,x,y)$ of a certain (finite!) degree in $Q$, that can be found automatically. This entails a certain non-linear recurrence
for the coefficients, and since we already
know that the initial conditions are all identically $0$, it is $0$ for ever after. A more efficient way is to use
the {\it holonomic ansatz} [Z] [K]. Eq. $(3)$ implies that its (unique) solution, that we also  call $\psi(x,y)$ (and we want to prove that
it is identical to the unique solution of $(1)$) is {\it holonomic}, i.e. satisfies linear differential equations with polynomial coefficients in {\bf both} $x$, and $y$.
This is also true for $g(x)=\psi(x,0)$. Using the `holonomic calculator' [K], we can get a holonomic description of
$Q(\psi(x,y),\psi(x,0),x,y)$ and checking sufficiently many initial values (coefficients in $x$) and verifying that they are $0$, it follows
by induction that all the coefficients are $0$ for ever after, i.e  $(3)$ implies $(1)$.

Now that we have a rigorous proof of $(3)$, plug-in $y=0$, and you get a rigorous proof of the previously-only-conjectured eq. $(1)$.

Now using Maple's {\tt gfun[algeqtodiffeq]} followed by {\tt gfun[diffeqtorec]} (or using our own home-made procedures  {\tt algtorec(F,P,x,n,N)}), we get
a linear recurrence equation with polynomial coefficients for the desired sequence of coefficients of $g(x)=\psi(x,0)$, namely the sequence $\{\psi_{n,0}\}$.

Tutte was also interested in the higher coefficients, in $y$, of the formal power series $\psi(x,y)$, and we can do it as well.

{\bf Sample Output files}

To see several examples of (rigorously proved) recurrences for the Taylor coefficients of $\psi(x,0)$ where $\psi(x,y)$ satisfies
other such functional equations, see the front of this article.

{\tt https://sites.math.rutgers.edu/\~{}zeilberg/mamarim/mamarimhtml/tutte.html} \quad .

Readers are more than welcome to experiment with their own favorite functional equation. The  function call is

{\tt Paper0(FE,psi,g,x,y,MaxC,n,G)} \quad,

where

$\bullet$ {\tt FE} is a polynomial with integer (or symbolic!, but beware things get very slow then)  coefficients in the variables {\tt psi, g, x,y}
(where {\tt g} is short for $\psi(x,0)$).

$\bullet$ {\tt psi,g,x,y} are symbols.

$\bullet$ {\tt MaxC} is a positive integer indicating the maximal complexity (order+degree) of the desired
recurrence for $\psi_{n,0}$ you are willing to take. You can always make it larger.

$\bullet$  {\tt G} is a large positive integer.

The {\bf output } is a computer-generated paper with 

$\bullet$ The algebraic equation satisfied by $g(x)=\psi(x,0)$ (that is initially guessed, but then rigorously proved {\it a posteriori}).

$\bullet$ The algebraic equation satisfied by $\psi(x,y)$  .

$\bullet$ The linear recurrence equation with polynomial coefficients satisfied by the coefficients of $g(x)$, alias the sequence $\{\psi_{n,0}\}$

(if there is no such recurrence with degree+order $\leq$ {\tt MaxC}, then, of course, it is not displayed. You are always welcome to
make {\tt MaxC} larger.

$\bullet$ If there is such a recurrence, then it gives the exact value of $\psi_{G,0}$ .

Enjoy!

{\bf References}

[GZ] Ira Gessel and Doron Zeilberger, {\it An Empirical Method for Solving (rigorously!) Algebraic Functional Equations Of the Form F(P(x,t), P(x,1),x,t)=0},
Personal Journal of Shalosh B. Ekhad and Doron Zeilberger, Dec. 28, 2014. \hfill\break
{\tt https://sites.math.rutgers.edu/\~{}zeilberg/mamarim/mamarimhtml/funeq.html} \quad .  \hfill\break
arxiv: {\tt https://arxiv.org/abs/1412.8360} \quad .

[K] Christoph Koutschan, {\it Advanced applications of the holonomic systems approach}, PhD thesis, \hfill\break
{\tt http://www.koutschan.de/publ/Koutschan09/thesisKoutschan.pdf} \quad .

[SZ]  Bruno Salvy and Paul Zimmermann, {\it GFUN: a Maple package for the manipulation of generating and holonomic functions in one variable},
ACM Transactions on Mathematical Software {\bf 20} (1994), 163-177. \hfill\break
{\tt https://dl.acm.org/doi/pdf/10.1145/178365.178368} \quad.

[T] W. T. Tutte, {\it A census of planar triangulations}, Canad. J. Math. {\bf 14} (1962), 21-38. \hfill\break
{\tt https://sites.math.rutgers.edu/\~{}zeilberg/akherim/tutte1962a.pdf} \quad .

[Z] Doron Zeilberger, {\it A Holonomic Systems approach to special functions}, J. Computational and Applied Math {\bf 32} (1990), 321-368. \hfill\break
{\tt https://sites.math.rutgers.edu/\~{}zeilberg/mamarim/mamarimhtml/holonomic.html} \quad .

\bigskip
\hrule
\bigskip
Shalosh B. Ekhad, c/o D. Zeilberger, Department of Mathematics, Rutgers University (New Brunswick), Hill Center-Busch Campus, 110 Frelinghuysen
Rd., Piscataway, NJ 08854-8019, USA. \hfill\break
Email: {\tt ShaloshBEkhad at gmail dot com}   \quad .
\bigskip
Doron Zeilberger, Department of Mathematics, Rutgers University (New Brunswick), Hill Center-Busch Campus, 110 Frelinghuysen
Rd., Piscataway, NJ 08854-8019, USA. \hfill\break
Email: {\tt DoronZeil at gmail  dot com}   \quad .
\bigskip
{\bf Exclusively published in the Personal Journal of Shalosh B. Ekhad and Doron Zeilberger, and arxiv.org}
\bigskip
{\bf March 8, 2024} \quad .

\end